\documentclass[letterpaper, 10 pt, journal, twoside]{IEEEtran}

\usepackage{amsmath}
\usepackage{placeins} 
\usepackage{breqn}
\usepackage{mathtools}
\usepackage{bm}
\usepackage{color, soul}
\usepackage{tabls}
\usepackage{units}
\usepackage{nicefrac}
\usepackage{cite}
\usepackage{tikz}
\usetikzlibrary{calc}
\usepackage[european]{circuitikz}
\usepackage{hyperref}
\usepackage{placeins}
\newcommand{\subparagraph}{}
\usepackage{nicefrac}
\usepackage{booktabs}
\usepackage{tabu}

\usepackage{xcolor}

\usepackage{titlesec}
\titlespacing*{\section}{0pt}{8pt}{3pt}
\titlespacing*{\subsection}{0pt}{3pt}{3pt}

\newcommand{\set}[1]{\mathcal{#1}} 
\newcommand{\Cov}{\mathrm{Cov}} 
\newcommand{\Var}{\mathrm{Var}} 
\newcommand{\epssum}{\tilde{\bm{\epsilon}}}

\usepackage{amsmath,amsfonts,amsthm}

\begin{document}

\pagestyle{empty}
\bstctlcite{IEEE:BSTcontrol} 

\title{Data-Driven Distributionally Robust Optimal Power Flow for Distribution Systems
}

\author{
    Robert Mieth, \emph{Student Member, IEEE}, 
    Yury Dvorkin, \emph{Member, IEEE}
    \vspace{-0.5ex}
}


\maketitle
\thispagestyle{empty}

\begin{abstract}
Increasing penetration of distributed energy resources complicate operations of electric power distribution systems by amplifying volatility of nodal power injections. On the other hand, these resources can  provide additional control means to the distribution system operator (DSO). This paper takes the DSO perspective and leverages a data-driven distributionally robust decision-making framework to overcome the uncertainty of  these injections and its impact on the distribution system operations. 
We develop an AC OPF formulation for radial distribution systems based on the  \emph{LinDistFlow} AC power flow approximation and exploit \emph{distributionally robust optimization} to immunize the optimized decisions against uncertainty in the probabilistic models of forecast errors obtained from the available observations. The model is reformulated to be computationally tractable and tested on multiple IEEE distribution test systems. We also release the code supplement that implements the proposed model in Julia and can be used to reproduce our numerical results, see [27].
\end{abstract}

{}

\IEEEpeerreviewmaketitle

\section{Introduction}

\IEEEPARstart{H}{istorically}, power transmission and distributions systems have been operated separately by the transmission system operator (TSO) and distribution system operator (DSO). Although both systems are interdependent and can be coordinated \cite{Bragin_2018pes}, the operational priorities and thus challenges vary. The TSO aims to continuously maintain  nodal power balances and mainly struggles with avoiding transmission overloads. On the other hand, the DSO is more focused on complying with nodal voltage limits in distribution systems and on following the pre-defined power exchanges with the TSO. Both operational paradigms are under pressure handling the constantly increasing volatility of power generation due to rising numbers of distributed energy resources (DERs) and aging infrastructure \cite{1338121}. While the integration of DERs is a policy priority in many jurisdictions, reliability and safety concerns may limit the techno-economic benefits of these resources. This paper takes the DSO perspective and aims to facilitate further integration of DERs by leveraging a data-driven distributionally robust decision-making framework to overcome the impacts of uncertain power injections on distribution systems.

Optimal Power Flow (OPF) tools  are routinely used to schedule and continuously dispatch controllable generators and loads to balance the system with minimal costs and losses with respect to the systems technical constraints (e.g. limits on generation outputs, voltage magnitudes and line flows). The inability to meet these limits may cause significant violations and thus may lead to voltage instability and, eventually, to cascading failures \cite{Bienstock_Chance_2014}. To avoid violating these limits in the presence of uncertain power injections, \emph{stochastic programming} and especially \emph{chance constrained optimization} have been leveraged for uncertainty-aware OPF models. The majority of such models have been designed for transmission systems, e.g. \cite{Bienstock_Chance_2014,Roald_Corrective_2017,Roald_Risk_2014}, and thus are tailored towards their operating needs. The studies in \cite{Bienstock_Chance_2014,Roald_Corrective_2017,Roald_Risk_2014} present  a risk-controlled chance-constrained OPF model (CCOPF) leveraging the DC-load-flow (DCLF) linearization, a common power flow approximation for originally non-convex and nonlinear AC power flow equations. These DC-CCOPF models are proven to effectively trade-off the likelihood of constraint violations%
and the security cost to avoid these violations. 
On the other hand, DERs are primarily located in distribution systems, where they mainly complicate voltage regulation, \cite{Hassan_Chance_2017}, and therefore the DCLF approximation is not technically suitable since it parametrizes voltage magnitudes at rated values. Dall'Anesse et al. \cite{DallAnese_Chance_2017} propose an AC-CCOPF model for distribution systems that linearizes AC power flows to enforce compliance with nodal voltage limits. 

The models in \cite{Bienstock_Chance_2014,Roald_Corrective_2017,Roald_Risk_2014,Hassan_Chance_2017,DallAnese_Chance_2017} exploit the common assumption of modeling uncertain nodal power injections using standard probability distributions with \textit{known parameters} such as mean and variance. This assumption does not hold in practice (e.g. see a data-driven wind power study in \cite{7268773}). The underlying distribution is never observable, but must be inferred from data. However, a stochastic program tuned towards a given data set often performs poorly when confronted with a different data set, even if it is drawn from the same distribution~\cite{Esfahani_Data_2015}. Instead of immunizing optimal solutions against worst-case observations that are available (data-robust methods) \emph{distributionally robust optimization} takes the worst-case over a family of distributions that are supported by the sample data \cite{Esfahani_Data_2015,Blanchet_Doubly_2017,Delage_Distributionally_2010,an_Ambiguous_2006}. In the OPF context, distributionally robust optimization has been extensively studied in the context of the transmission OPF models, e.g. \cite{Duan_Distributionally_2017,Lubin_A_2016,Guo_Stochastic_2017}. These studies use the DCLF approximation to represent power flows in a form suitable for introducing chance constraints and assume knowledge of first and second order distribution moments \cite{Roald_Security_2015}.
 
This paper presents a data-driven distributionally robust AC CCOPF (DR-AC-CCOPF) for distribution systems,  building on the results in \cite{Duan_Distributionally_2017,Lubin_A_2016,Guo_Stochastic_2017} for transmission systems. The uncertainty on nodal power injections is caused by net load fluctuations, which are defined in this work as the difference between the demand and power output of uncontrollable DERs (``behind-the-meter systems") at each node. To represent AC power flows, we leverage the established \emph{LinDistFlow} \cite{Baran_Optimal_1989}, a linear representation of the AC power flows that neglects power losses, but accounts for apparent flows and voltage magnitudes. The chance constraints and the objective are then recast as second-order expressions, thus leading to a convex (quadratic) OPF program. We then robustify the formulation by introducing an ambiguous probability distribution of the uncertain input via a distributional uncertainty set (cf. \cite{Esfahani_Data_2015}). Our case study corroborates the usefulness and scalability of the proposed DR-AC-CCOPF.

\section{AC-CCOPF in Radial Distribution Networks}
\label{sec:cc_radial_opf}

\subsection{Preliminaries}

We consider a radial distribution system represented by graph $\Gamma (\set{N}, \set{E})$, where $\set{N}$ and $\set{E}$ are the sets of nodes (buses) and edges (lines), respectively. The graph is a tree with the root indexed as $0$ and $\set{N}^+ \coloneqq \set{N} \setminus \{0\}$ is the set of all non-root nodes. The root node is usually referred to as the \emph{substation}, as it connects the distribution system with the transmission system. Each node is characterized by its net demand, its active and reactive power ($d_i^P$ and $d_i^Q$, $i \in \set{N}$) and its voltage magnitude $v_i \in [v_i^{min}, v_i^{max}], i \in \set{N}$, where $v_i^{min}$ and $v_i^{max}$ are the upper and lower voltage limits. Also, we introduce $u_i = v_i^2, i \in \set{N}$. If a node is equipped with a controllable DER (e.g. cogeneration resource), we model active and reactive power generation $g_i^P \in [g_i^{P,min}, g_i^{P,max}]$ and $g_i^Q \in [g_i^{Q,min},g_i^{Q,max}]$, $i \in \set{G} \subseteq \set{N}$, where $g_i^{P,max}$ and $g_i^{Q,max}$ are the output limits. The same logic is applied for the substation which will be treated as a generator that provides power by buying it from the transmission grid.
Each node is associated with an ancestor (or parent) node $\set{A}_i$ and a set of children nodes $\set{C}_i$. Since $\Gamma $ is radial, it is \mbox{$|\set{A}_i| = 1, \forall i \in \set{N}^+$} and all edges $i \in \set{E}$ are indexed by $\set{N}^+$.

\begin{figure}[!b]
\centering
\includegraphics[width=0.7\linewidth]{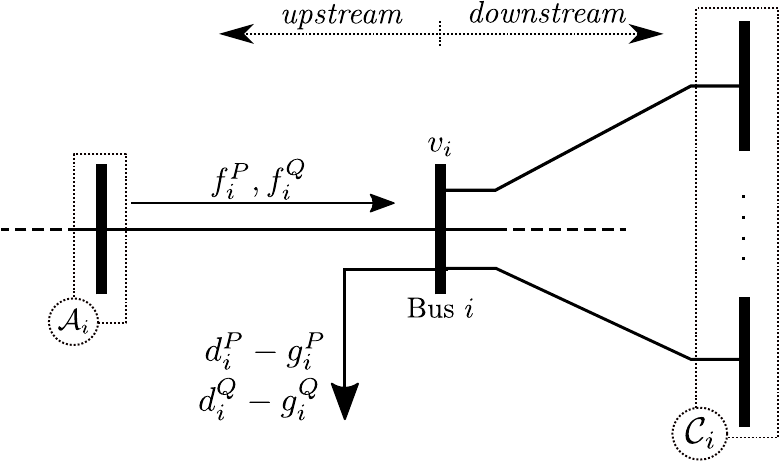}
\caption{Main notations of the radial distribution system.}
\label{fig:notation_illustration}
\end{figure}

The active and reactive power flows are denoted as $f_i^P$ and $f_i^Q$, $i \in \set{E}$, where $i$ is the index of the \emph{downstream} node of edge~$i$, i.e. the node at the receiving end of edge~$i$. Accordingly, the sending node of edge $i$ is called called \emph{upstream}. Each edge $i$ has resistance $R_i$, reactance $X_i$ and apparent power flow limit $S_i^{max}, i \in \set{E}$.  
Figure \ref{fig:notation_illustration} illustrates the introduced notations. 
In the following, we leverage the \emph{LinDistFlow} \cite{Turitsyn_Local_2010,Baran_Optimal_1989} approximation as:
\begin{align}
        (d_i^p - g_i^p) + \sum_{j \in \set{C}_i} f_j^p = f_{i}^p, \quad p \in \{P,Q\} 
            &\quad \forall i \in \set{N}^+ \label{eq:first_lindist_const}\\
        u_{\set{A}_i} - 2(f_i^P R_i + f_i^Q X_i) = u_i, 
            &\quad \forall i \in \set{N}^+ \label{eq:lindist_volt}
\end{align}
\noindent
where $u_0$ is assumed to be the base voltage on the secondary side of the substation, i.e. $u_0 = \unit[1]{\text{p.u.}}$. The distribution network and controllable resources are constrained as: 
\begin{align}
    0 \leq g_i^P \leq g_i^{P,max},     
        &\quad \forall i \in \set{G} \label{eq:lindist_gP_const}\\
    -g_i^{Q,max} \leq g_i^Q \leq g_i^{Q,max},    
        &\quad \forall i \in \set{G} \\
    u_i^{min} \leq u_i \leq u_i^{max},  
        &\quad \forall i \in \set{N} \label{eq:last_lindist_const} \\
    (f_i^P)^2 + (f_i^Q)^2 \leq (S_i^{max})^2,    
        &\quad \forall i \in \set{E}.
\end{align}

\subsection{Chance Constraint OPF for LinDistFlow}
An OPF problem that takes the DSO perspective seeks to dispatch controllable generators in the least-cost manner subject to the system's physical constraints. However, the nodal net load injections in \eqref{eq:first_lindist_const} are unknown ahead of time. We therefore define the uncertain net load injections as $\bm{d}_i \coloneqq \bar{d}_i + \bm{\epsilon}_i$, where $\bar{d}_i$ is a given forecast and $\bm{\epsilon}_i$ is a normally distributed forecast error with mean $\mathbb{E}(\bm{\epsilon}_i) = \mu_i$, standard-deviation $\sigma_i$ and variance $\Var(\bm{\epsilon}_i) = \sigma_i^2$. By assuming that the forecast is not systematically flawed we get $\mu_i = 0$ and thus $\mathbb{E}(\bm{d}_i) = \bar{d}_i$. Accordingly, for non-zero mean one could use $\mathbb{E}(\bm{d}_i) = \bar{d}_i + \mu_i$. However, following the relevant literature e.g. \cite{Bienstock_Chance_2014,Lubin_A_2016} we will not consider this case.
Cf. \cite{7268773, Roald_Security_2015} for a comprehensive discussion on the distribution of net load forecast errors. We also assume independent forecast errors at different buses, following \cite{Bienstock_Chance_2014}. Note that this assumption is in line with the conservatism of the distributionally robust optimization.
We will use bold symbols for uncertain parameters. 

For the distinction between forecast errors in active and reactive power components of the nodal net load injections, there are two reasonable treatment methods:
 
\begin{enumerate}
    \item The uncertain nodal injections result from an unpredictable combination of various appliances leading to an equally unpredictable power factor at the bus. If the modeled variance in active and reactive power demand is sufficiently small so that the realized power factor is physically consistent, the random errors can be treated independently as $\bm{\epsilon}_i^P$ and $\bm{\epsilon}_i^Q$. 
    \item The uncertain nodal injections result from the unpredictable utilization of large appliances with a constant power factor. In this case, $\bm{\epsilon}_i^P = \bm{\epsilon}_i$ and $\bm{\epsilon}_i^Q = \bm{\epsilon}_i \tan \varphi_i$ with \mbox{$\cos \varphi_i = \nicefrac{\bar{d}^P_i}{\sqrt{(\bar{d}^P_i)^2 + (\bar{d}^Q_i)^2}}$} being the power factor.  
\end{enumerate}
For the sake of generality, we use $\bm{\epsilon}_i^P$ and $\bm{\epsilon}_i^Q$ in the following notations without assuming a specific treatment method.

When the uncertainty materializes, the controllable generators must compensate for net load deviations from the forecast values. Following the established policy of affine real-time control \cite{Bienstock_Chance_2014,on_and_Probabilistic_1974,Roald_Corrective_2017}, the distribution of the load deviation among the controllable generators is decided by variable $\alpha_i \geq 0, i \in \set{G}$, such that:
\begin{align}
    \bm{g}_i^P = \bar{g}_i^P + \alpha_i \epssum^P, 
        &\quad i \in \set{G} \label{eq:uncertain_gP}\\
    \bm{g}_i^Q = \bar{g}_i^Q + \alpha_i \epssum^Q, 
        &\quad i \in \set{G}. \label{eq:uncertain_gQ}
\end{align}
\noindent
where $\epssum^P = \sum_{i \in \set{N}} \bm{\epsilon}^P_i$ is the sum of forecast errors ($\epssum^Q$ by analogy) and $\bar{g}_i^P$, $\bar{g}_i^Q$ are the dispatch decisions of controllable generators based on the forecast. Note that $\alpha_i$ is a decision variable and that $\sum_{i \in \set{G}} \alpha_i = 1$. We additionally formalize the vector $\alpha$ as:
\begin{equation*}
    \alpha_{(i)} = \begin{cases}
        \alpha_i, &\parbox[t]{.2\textwidth}{if $i \in \set{G}$} \\
        0,        & \text{otherwise} 
    \end{cases} \quad i \in \set{N}.
\end{equation*}

Similarly to \eqref{eq:uncertain_gP} and \eqref{eq:uncertain_gQ}, we now derive expressions for the edge power flows and nodal voltages affected by the uncertain nodal injections. We introduce matrix $A: l \times b$, where $l \coloneqq |\set{E}|$ is the number edges and $b \coloneqq |\set{N^+}|$ is the number nodes. We leverage the concept of \emph{power transmission distribution factors} (PTDF) \cite{Christie_Transmission_2000} to map the change of load at every node to the change of edge power flow. Note that since the distribution system is radial, each edge has to carry the complete net load of all its downstream nodes. Hence, is $A$ constructed as:
\begin{equation*}
    a_{ij} \coloneqq \begin{cases}
        1, & \parbox[t]{.2\textwidth}{if line $i$ is part of the path from root to bus $j$} \\
        0, & \text{otherwise}
    \end{cases} \quad i \in \set{E},~j \in \set{N^+},  
\end{equation*} 
\noindent
where $a_{ij}$ denotes the entry in column $j$ of row $i$ of the \mbox{matrix $A$}. 

Using $A$ and $\alpha$ as defined above, the uncertain edge power flow can be written as: 
\begin{align}
    \bm{f}_i^P = \bar{f}_i^P + a_{i*} (\bm{\epsilon}^P - \alpha \bm{\tilde{\epsilon}}^P), \label{eq:uncertain_fP}\\
    \bm{f}_i^Q = \bar{f}_i^Q + a_{i*} (\bm{\epsilon}^Q - \alpha \bm{\tilde{\epsilon}}^Q), \label{eq:uncertain_fQ}  
\end{align}
\noindent
where $\bm{\epsilon}^{P}: b \times 1 \coloneqq (\bm{\epsilon}^P_1,...,\bm{\epsilon}_n^P)$ is the vector of active power forecast errors of all nodes ($\bm{\epsilon}^Q$ by analogy) and $\bar{f}_i^P$ and $\bar{f}_i^Q$ are the flows for the forecasted nodal net load injections. Row-vector $a_{i*}: 1 \times b$ corresponds to the $i$-th row of $A$. 
Using \eqref{eq:lindist_volt} with respect to \eqref{eq:uncertain_fP} and \eqref{eq:uncertain_fQ}, we find the uncertain voltage magnitudes squared  as:
\begin{align}
    \bm{u}_i &= \bm{u}_{\set{A}_i} - 2(R_i\bm{f}_i^P + X_i\bm{f}_i^Q)  \label{eq:uncertain_voltage} \\
             &= \bar{u}_i - 2 a_{*i}^\top (R \circ A(\bm{\epsilon}^P - \alpha \bm{\tilde{\epsilon}}^P) \nonumber 
             + X \circ A(\bm{\epsilon}^Q - \alpha \bm{\tilde{\epsilon}}^Q)), 
\end{align} 
\noindent
where $(\circ)$ denotes the Schur-product and $R$, $X$ are the vectors of line resitances and reactances respectively.
Given the uncertain nodal net load injections, we can now formulate the AC-CCOPF problem that accounts for the propagation of the uncertainty within a given radial distribution system and its impact on power flow and voltage limits: 
\begin{subequations}
\begin{align}
   \min_{(\bar{g}, \alpha, \bar{f}, \bar{u})} \mathbb{E}[ f(\bar{g}^P, \bar{g}^Q, \bm{\epsilon}^P, \bm{\epsilon}^Q)] \label{eq:_obj_function}\\
    \sum_{i \in \set{G}} \alpha_i = 1 \qquad &\\
    (\bar{d}_i^p - \bar{g}_i^p) + \sum_{j \in \set{C}_i} \bar{f}_j^p 
        = \bar{f}_i^P, \quad p \in \{P,Q\} \quad &\forall i \in \set{E}\\
    \bar{u}_{\set{A}_i} - 2(\bar{f}_i^P R_i + \bar{f}_i^Q X_i) 
        = \bar{u}_{i}, \quad &\forall i \in \set{N}^+ \\
    \mathbb{P}(g_i^{P,min} \leq \bm{g}_i^P \leq g_i^{P,max}) \geq (1-2\eta_g), \quad 
        &\forall i \in \set{G}_c, \label{eq:first_cc} \\
    \mathbb{P}(g_i^{Q,min} \leq \bm{g}_i^Q \leq g_i^{Q,max}) \geq (1-2\eta_g), \quad
        &\forall i \in \set{G}_c, \label{eq:reactive_gen_cc} \\
    \mathbb{P}(u_i^{min} \leq \bm{u}_i \leq u_i^{max}) \geq (1-2\eta_v), \quad
         &\forall i \in \set{N}. \label{eq:last_cc}\\
 \quad \eqref{eq:uncertain_gP} - \eqref{eq:uncertain_voltage} 
    \end{align}
    \label{eq:cc_ropf}
\end{subequations} 
\noindent where \eqref{eq:first_cc}-\eqref{eq:last_cc} are the chance constrained equivalents of \eqref{eq:lindist_gP_const}-\eqref{eq:last_lindist_const}, whereas parameters $\eta_g, \eta_v, \eta_f \in (0,\nicefrac{1}{2}]$ set the probability limit on constraint violation. Generally, $\eta_g, \eta_v, \eta_f$ are small and confidence levels, e.g. $(1-\eta_g)$, of \unit[95 - 99]{\%} are common. The factor $2$ is necessary to be consistent in using two-sided chance constraints as discussed in Section \ref{ssec:conic_reformulation}. The operator $\mathbb{E}(\cdot)$ denotes the expected value over the forecast error distribution and $\mathbb{P}(X)$ the probability of event $X$. 

\textit{Remark:} The AC-CCOPF (12) does not impose power flow limits on $\boldsymbol{f}_{t,l}^p$ and $\boldsymbol{f}_{t,l}^q$, because real-life distribution systems are typically voltage-constrained and power flow limits can be disregarded (e.g. future distribution operations are expected to have ``\textit{bounds on system frequency, voltage levels, and DER capacities}'' \cite{taft2015reference}. In \cite{lubin2018acccopf}, we propose an approach to enforce chance-constrained apparent power limits by a suitable inner approximation.

\subsection{Expected Cost}
\label{ssec:expected_cost}

Let $c_i(\cdot), i \in \set{G},$ be the quadratic cost function \cite{Bienstock_Chance_2014} of each controllable generator: 
\begin{equation}
     c_i(\bm{g}^P) = c_{i2} (\bm{g}^P)^2 + c_{i1} \bm{g}^P + c_{i0}. 
\end{equation}
\noindent 
Then, the system-wide generation cost is:
\begin{subequations}
\begin{align}
    f_P(\bar{g}^P, \bm{\epsilon}^P) = \sum_{i \in \set{G}} c_i(\bm{g}^P) = c_2 (\tilde{\bm{\epsilon}}^P)^2 + c_1 \tilde{\bm{\epsilon}}^P + c_0, 
\end{align}

\noindent
where
\begin{equation}
    \begin{cases}
        c_2 &= \sum_{i \in \set{G}} c_{i2} \alpha_i^2 \\
        c_1 &= \sum_{i \in \set{G}} \left(2 c_{i2} \alpha_i \bar{g}^P_i + c_{i1} \alpha_i \right) \\
        c_0 &= \sum_{i \in \set{G}} \left(c_{i2} (\bar{g}^P_i)^2 + c_{i1}\bar{g}^P_i + c_{i0} \right). 
    \end{cases}
\end{equation}
\end{subequations}

\noindent Under the  zero-mean assumption on the forecast error invoked above and the independence of the forecast errors at different nodes, i.e. $\Cov(\bm{\epsilon}_i, \bm{\epsilon}_j) = 0, \forall i \neq j \in \set{N}$, we get the following expression for the expected cost in \eqref{eq:_obj_function}:
\allowdisplaybreaks
\begin{align}
    &\mathbb{E}\left[f_P(\bar{g}^P, \bm{\epsilon}^P)\right] \nonumber = \mathbb{E}\left[ c_2 (\tilde{\bm{\epsilon}}^P)^2 + c_1 \tilde{\bm{\epsilon}}^P + c_0  \right] \nonumber \\
        &= c_2\mathbb{E}\left[(\tilde{\bm{\epsilon}}^P)^2 \right] + c_1 \mathbb{E}\left[\tilde{\bm{\epsilon}}^P \right] + c_0  \nonumber \\
        &= c_2 \Var(\tilde{\bm{\epsilon}}^P) + c_0 \nonumber \\
        &= \sum_{i\in \set{G}} \left[ c_{i2} \left(\alpha_i^2 \Var(\tilde{\bm{\epsilon}}^P) + (\bar{g}^P_i)^2 \right) + c_{i1}\bar{g}^P_i + c_{i0} \right],
        \label{eq:expected_costs}
\end{align}
\allowdisplaybreaks[0]
\noindent where  $\Var(\tilde{\bm{\epsilon}}^P) = \sum_{i \in \set{N}} (\sigma_i^P)^2$. 
Since it is common practice of real-life utilities (e.g. Con Edison of NY) to not price reactive power we will forgo to derive an expression for expected costs induced by the reactive power component. The process, however, is the same.  

\subsection{Conic Reformulation of the Chance Constraints}
\label{ssec:conic_reformulation}

Following \cite{Bienstock_Chance_2014,Roald_Security_2015}, the chance constraints in \eqref{eq:cc_ropf} can be reformulated into computationally tractable second-order conic constraints. 
We also refer interested readers to  \cite[Chapter~4.4.2]{Boyd_Convex_2004} for more details on conic programming.

For any normally distributed random variable \mbox{$X \sim \mathrm{Norm}(\mu, \sigma^2)$}, the inequality $\mathbb{P}(X \leq x^{max}) > \eta$ holds if and only if $x^{max} \geq \mu + z_\eta \sigma$, where $z_\eta \coloneqq \Phi^{-1}(1-\eta)$ is the $(1-\eta)$-quantile of a standard normal distribution %
\footnote{The same transformation is possible with different symmetric and unimodal distributions by changing $\Phi^{-1}(1-\eta)$ to a general inverse cumulative distribution function $f^{-1}(1-\eta)$, cf. \cite{Roald_Security_2015} for details.}. 
This transformation can be exploited for the upper-bound chance constraint in \eqref{eq:first_cc} and \eqref{eq:reactive_gen_cc} as follows:
\begin{subequations}
\begin{align}
    g_i^{p,max} &\geq \mathbb{E}(\bm{g}_i^{p}) + z_{\eta_g} \sqrt{\Var(\bm{g}_i^p)} \nonumber \\
            &= \bar{g}_i^{p} + z_{\eta_g} \alpha_i \sqrt{\sum_{j=1}^b \Var(\bm{\epsilon}_j^p)}, 
            \quad p \in \{P,Q\}. \label{eq:conic_gmax}
\end{align}
\noindent
where $\alpha_i \geq 0$ as introduced in \eqref{eq:uncertain_gP} and \eqref{eq:uncertain_gQ}.

Similarly, the lower-bound chance constraints of \eqref{eq:first_cc} and \eqref{eq:reactive_gen_cc} can be reformulated as:
\begin{equation}
     -g_i^{p,min} \geq -\bar{g}_i^{p} + z_{\eta_g} \alpha_i \sqrt{\sum_{j=1}^b \Var(\bm{\epsilon}_j^p)}, \quad p \in \{P,Q\}. \label{eq:conic_gmin}
 \end{equation} 
\label{eq:conc_g_limits}
\end{subequations}

Finally the voltage constraints can be recast as:
\begin{subequations}
\begin{align}
    u_i^{max} \geq& \bar{u}_i + z_{\eta_v} \sqrt{\Var(\bm{u}_i)} \\
    -u_i^{min} \geq& -\bar{u}_i + z_{\eta_v} \sqrt{\Var(\bm{u}_i)}
\end{align}
\begin{align}
    &\Var(\bm{u_i}) = 4 \sum_{k=1}^l a_{ki}  \left[ \left( R_k^2 \sum_{j=1}^b a_{kj} (\Var(\bm{\epsilon}_j^P) + \alpha_j^2 \Var(\tilde{\bm{\epsilon}}^P)) \right) \right. \nonumber \\
    & \qquad \left. + \left( X_k^2 \sum_{j=1}^b a_{kj} (\Var(\bm{\epsilon}_j^Q) + \alpha_j^2 \Var(\tilde{\bm{\epsilon}}^Q)) \right)\right], 
\end{align}
\label{eq:var_u}
\end{subequations}

\noindent
where $\eta_v, \eta_g < \nicefrac{1}{2}$ so that $z_{\eta_v}, z_{\eta_g}>0$ and \eqref{eq:conc_g_limits} - \eqref{eq:var_u} are conic and thus computationally tractable (cf. \cite{Bienstock_Chance_2014,Boyd_Convex_2004}).

\section{Data-Driven DR-AC-CCOPF}
\label{sec:data-driven_dr_cc}

So far the presented formulations solve the stochastic problem with the assumption of perfect knowledge of the underlying distribution of the random variable. However, this true distribution can not be known exactly, but is only informed by finite observations of previous realizations. The distribution that we use to define both the expected costs and the chance constraints therefore is ambiguous over the available data. To robustify the formulation robust against  uncertain distributions, we redefine the AC-CCOPF objective as follows: 
\begin{equation}
   \min_{x} \sup_{\pi \in \set{U}} \mathbb{E}_{\pi}[f(x, \bm{\epsilon})],
\end{equation}
\noindent
where the set $\set{U}$ is the set of all distributions that are supported by the available data within a predefined level of confidence and $\mathbb{E}_{\pi}$ is the expectation taken with respect to the distribution~$\pi$. The task is to minimize the \emph{worst-case expectation} based on the distributional uncertainty set that has been inferred by the data. The quality of the available historic data does not affect our proposed model inherently. However, the accuracy of forecasts and volatility of net loads depend on meteorological and circumstantial externalities. Since stochastic optimization is especially powerful for short-term (e.g. intra-day) system dispatch, it is reasonable to assume the availability of detrended data based on similar conditions (e.g. type of day, season, etc).

\subsection{Definition of the Uncertainty Set}

Let $\set{H}(\bm{\epsilon}_i) \coloneqq \{\hat{\epsilon}_{i,t}\}_{t \leq N}, i \in \set{N}$ be the set of $N$ observed realizations of the forecast error at bus $i$ for either active or reactive power forecast, where $\hat{\epsilon}_i = \hat{d}_i - \bar{d}_i$. We introduce the circumflex ($\hat{~}$) to mark all values that are based on empirical data. As the forecast error is zero-mean, normally distributed,  the distribution can formally be defined via  the \emph{sample variance} of the empirical data at each bus as:
\begin{equation}
    \hat{\sigma}_i^2 = \frac{1}{N} \sum_{t \leq N}\hat{\epsilon}_{i,t}^2, \qquad i \in \set{N}.  \label{eq:sample_variance}
\end{equation}
Note that the zero mean assumption allows the calculation of $\hat{\sigma}^2_i$ with the full degree of freedom $N$ (as opposed to $N-1$ in the case of estimated mean). 

Although the sample variance is the \emph{minimum-variance unbiased estimator} (MVUE) of the unknown distribution, it will never resemble the true variance perfectly while $N<\infty$. It has been shown that the sample variance itself is a random-variable following a Chi-Square ($\chi^2$) distribution parameterized by the number of available samples $N$%
 (for details cf. e.g. \cite{pestman1998mathematical})
. We can use this property to define the set $\set{U}_{\sigma_i^2}$ based on $\hat{\sigma}_{i}$ as an interval that contains the true variance with probability $1-\xi$:
\begin{subequations}
\begin{equation}
   \set{U}_{\hat{\sigma}_{i}^2} = \left[ \hat{\zeta}_{i,l}, \hat{\zeta}_{i,h} \right],
\end{equation}
\noindent
where
\begin{align}
    \hat{\zeta}_{i,l} \coloneqq \frac{N\hat{\sigma}^2_{i}}{\chi^2_{N,(1-\xi)/2}}, \quad
    \hat{\zeta}_{i,h} \coloneqq \frac{N\hat{\sigma}^2_{i}}{\chi^2_{N,\xi/2}} \label{eq:sigma_upper_bound}.
\end{align}
\label{eq:sigma_interval}
\end{subequations}

It is $\chi^2_{N,\xi}$ the $\xi$-quantile of the $\chi^2$-distribution with $N$ degrees of freedom. Note that the width of the interval will decrease with the amount of available samples and that the $\chi^2$-distribution is not symmetric; the interval will therefore not be centered around $\hat{\sigma}^2$ (cf. Fig. \ref{fig:distribution_visualization} (c)).

\subsection{Worst-Case Expectation}

The expected system cost in \eqref{eq:expected_costs} is the sum of convex quadratic cost functions of individual generators, which includes the sum of the variances of the forecast errors. With respect to \eqref{eq:sigma_interval} the worst-case expectation is given as:
\begin{equation}
    \begin{aligned}
        &\sup_{\sigma^2_i \in \set{U}_{\hat{\sigma}_{i}^2}, i \in \set{N}} \mathbb{E}\left[f_P(\bar{g}^P, \bm{\epsilon}^P)\right] \\
        &= \sup_{\sigma^2_i \in \set{U}_{\hat{\sigma}_{i}^2}, i \in \set{N}} \sum_{i\in \set{G}} \left[ c_{i2} \left(\alpha_i^2 \Var(\tilde{\bm{\epsilon}}^P) + (\bar{g}^P_i)^2 \right) + c_{i0} \right] \\
        &= \sum_{i\in \set{G}} \left[ c_{i2} \left(\alpha_i^2 
        \sum_{i \in \set{N}} \hat{\zeta}_{i,h} 
        + (\bar{g}^P_i)^2 \right) + c_{i1}\bar{g}^P_i + c_{i0} \right].
    \end{aligned}
    \label{eq:wc_expected_cost}
\end{equation}

The linear relation between the sum of the individual error variances at the nodes and the expected cost leads to the upper bound of the uncertainty region $\hat{\zeta}_{i,h}$ as the worst case expectation of the objective function (cf.~\cite{Ben-Tal_Robust_2009}).

\subsection{DR-AC-CCOPF Formulation}

We  can now reformulate the AC-CCOPF problem such that the risk of constraint violations in the presence of uncertain load is minimized and also to account for our data-driven, incomplete knowledge of the underlying error distribution: 

\allowdisplaybreaks
\begin{subequations}
\begin{align}
&  \min_{(\bar{g}, \alpha, \bar{f}, \bar{u})} 
      \sum_{i\in \set{G}} \left[ c_{i2}^P \left(\alpha_i^2 
        \sum_{i \in \set{N}} \hat{\zeta}_{i,h}^P 
        \!+ \!(\bar{g}^P_i)^2 \right)\! +\! c_{i1}\bar{g}^P_i \!+\! c_{i0}^P \right] \\
& \sum_{i \in \set{G}} \alpha_i = 1 \\        
& (\bar{d}_i^P - \bar{g}_i^P) + \sum_{j \in \set{C}_i} \bar{f}_j^P = \bar{f}_i^P\quad \forall i \in \set{E}\\
& (\bar{d}_i^Q - \bar{g}_i^Q) + \sum_{j \in \set{C}_i} \bar{f}_j^Q = \bar{f}_i^Q\quad \forall i \in \set{E}\\
& \bar{u}_{\set{A}_i} - 2(\bar{f}_i^P R_i + \bar{f}_i^Q X_i) = \bar{u}_{i}\quad \forall i \in \set{N}^+ 
\\
& g_i^{p,max} \geq \bar{g}_i^{p} + z_{\eta_g} \alpha_i \sqrt{\sum_{k=1}^b \hat{\zeta}_{i,h}^p} \quad \forall i \in \set{G}, p \in \{P,Q\} \\
& -g_i^{p,min} \geq -\bar{g}_i^{p} + z_{\eta_g} \alpha_i \sqrt{\sum_{k=1}^b \hat{\zeta}_{i,h}^p} \quad \forall i \in \set{G}, p \in \{P,Q\} \\
& u_i^{max} \geq \bar{u}_i + z_{\eta_v} 2\sqrt{h_i(\alpha)} \quad \forall i \in \set{N} \\
& -u_i^{min} \geq -\bar{u}_i + z_{\eta_v} 2\sqrt{h_i(\alpha)} \quad \forall i \in \set{N}  \\
&h_i(\alpha) =  \sum_{k=1}^l a_{ki} \left[ \left( R_k^2 \sum_{j=1}^b a_{kj} (\hat{\zeta}_{j,h}^P + \alpha_j^2 \sum_{m=1}^b\hat{\zeta}_{m,h}^P))         \right) \right. \nonumber \\
    & \qquad \left. + \left( X_k^2 \sum_{j=1}^b a_{kj} (\hat{\zeta}_{j,h}^Q + \alpha_j^2 \sum_{m=1}^b\hat{\zeta}_{m,h}^Q)) \right)\right] \label{eq:dddr_vmax}  
 \end{align}
\end{subequations}
\allowdisplaybreaks[0]

As discussed in Section \ref{ssec:expected_cost}. the reactive power cost is neglected in this formulation. The resulting problem is quadratic due to $\alpha_i^2$ but solvable by  off-the-shelf solvers. 

\section{Case Study}
\label{sec:case_study}

\begin{figure}
\centering
\includegraphics[width=0.9\linewidth]{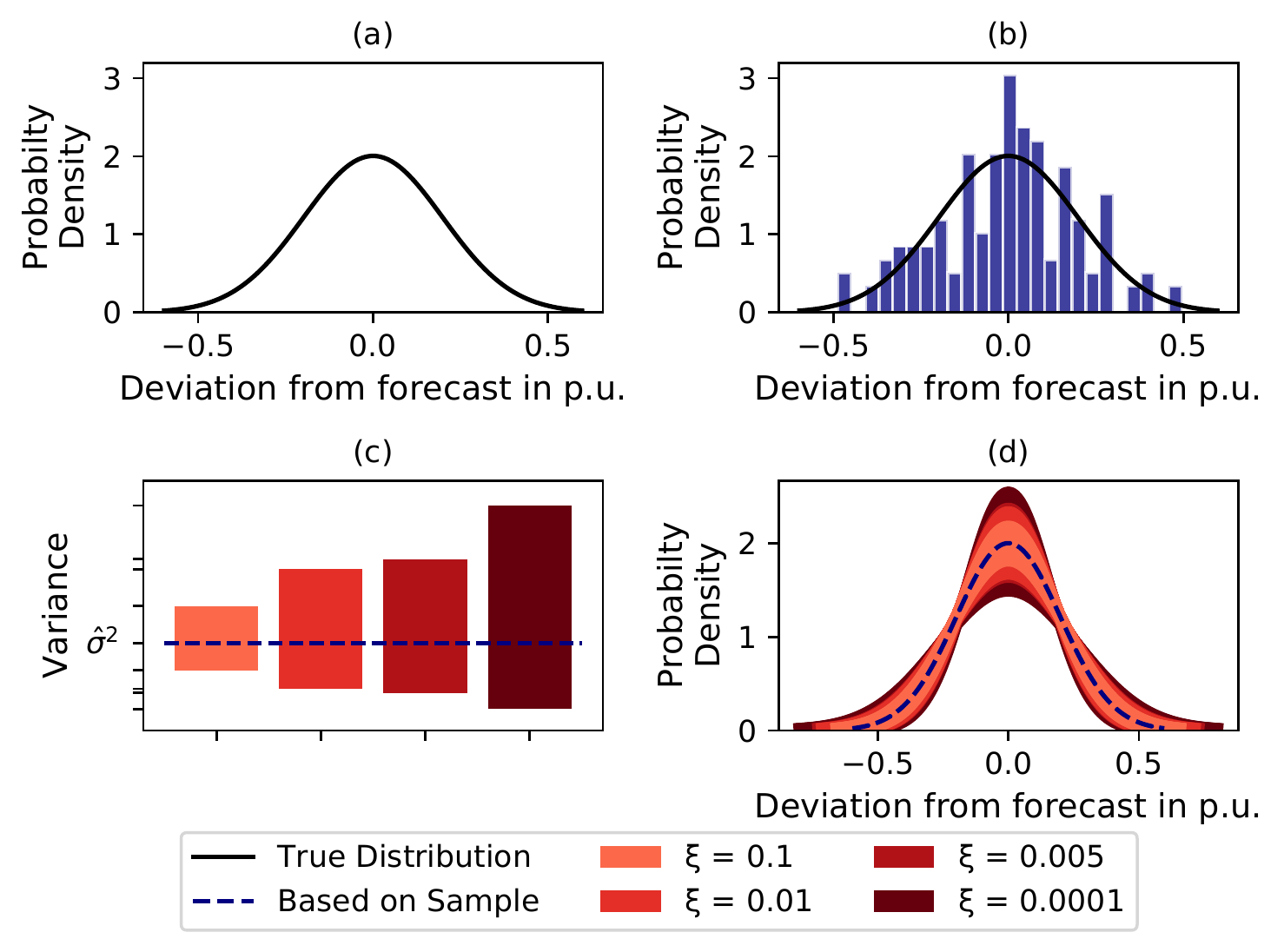}
\caption{Exemplary true distribution, historic observations and uncertainty sets for a bus with a load of \unit[1]{p.u.}: \mbox{\textbf{(a)} True} (unknown) error-distribution with a standard deviation of \unit[20]{\%} of the load. \mbox{\textbf{(b)} $N=100$ samples} drawn from the true error-distribution (blue bar-histogram). \textbf{(c)} Uncertainty intervals around sample variance $\hat{\sigma}_N$ for different $(1-\xi)$. \textbf{(d)} Set of possible distributions based on $\hat{\sigma}_i$ and the uncertainty intervals, respectively.}
\label{fig:distribution_visualization}
\end{figure}

We use the \mbox{15-bus} radial distribution system from \cite{Papavasiliou_Analysis_2017} and add two fully controllable generators at nodes 6 and 11 with the production cost of \$10/MWh each in addition to the substation at the root node, which supplies at \$50/MWh. These costs are selected to incentivize the use of distributed generators.  We assume that the net load forecasts are given for each node with the zero-mean forecast error and the standard deviation of $\sigma_i^P = 0.2 \bar{d}_i^P$,%
\footnote{This relation has been shown as a feasible assumption based on empirical data in \cite{Roald_Security_2015} and \cite{7268773}. The latter reference also shows how to overcome data misfits inflicted by assuming a normal distribution.}
as in Figure \ref{fig:distribution_visualization}(a).We use this true distribution to obtain $N=100$ error samples that maintain the node's power factor shown in the histogram of Figure \ref{fig:distribution_visualization}(b). In turn, we use these samples and \eqref{eq:sample_variance}-\eqref{eq:sigma_interval} to derive uncertainty intervals  for different $\xi$ as shown in Figure \ref{fig:distribution_visualization}(c). Figure \ref{fig:distribution_visualization}(d) displays the resulting distribution used to solve the AC-CCOPF (dashed line) and the uncertainty sets used to solve the DR-AC-CCOPF (colored areas). The AC-CCOPF and DR-AC-CCOPF formulations are then compared in terms of their constraint feasibility and operating cost. To this end, we solve each model and then test their solutions against 750 random samples generated from the true distribution in Fig.~\ref{fig:distribution_visualization}(a), which is sufficient to obtain stable empirical distributions in both cases. We implement the case study using the \emph{Julia JuMP} package~\cite{jump} and our code and input data can be downloaded from~\cite{cc_code}.

\subsection{In-Sample Evaluation}

\begin{figure}
\centering
\includegraphics[width=0.8\linewidth]{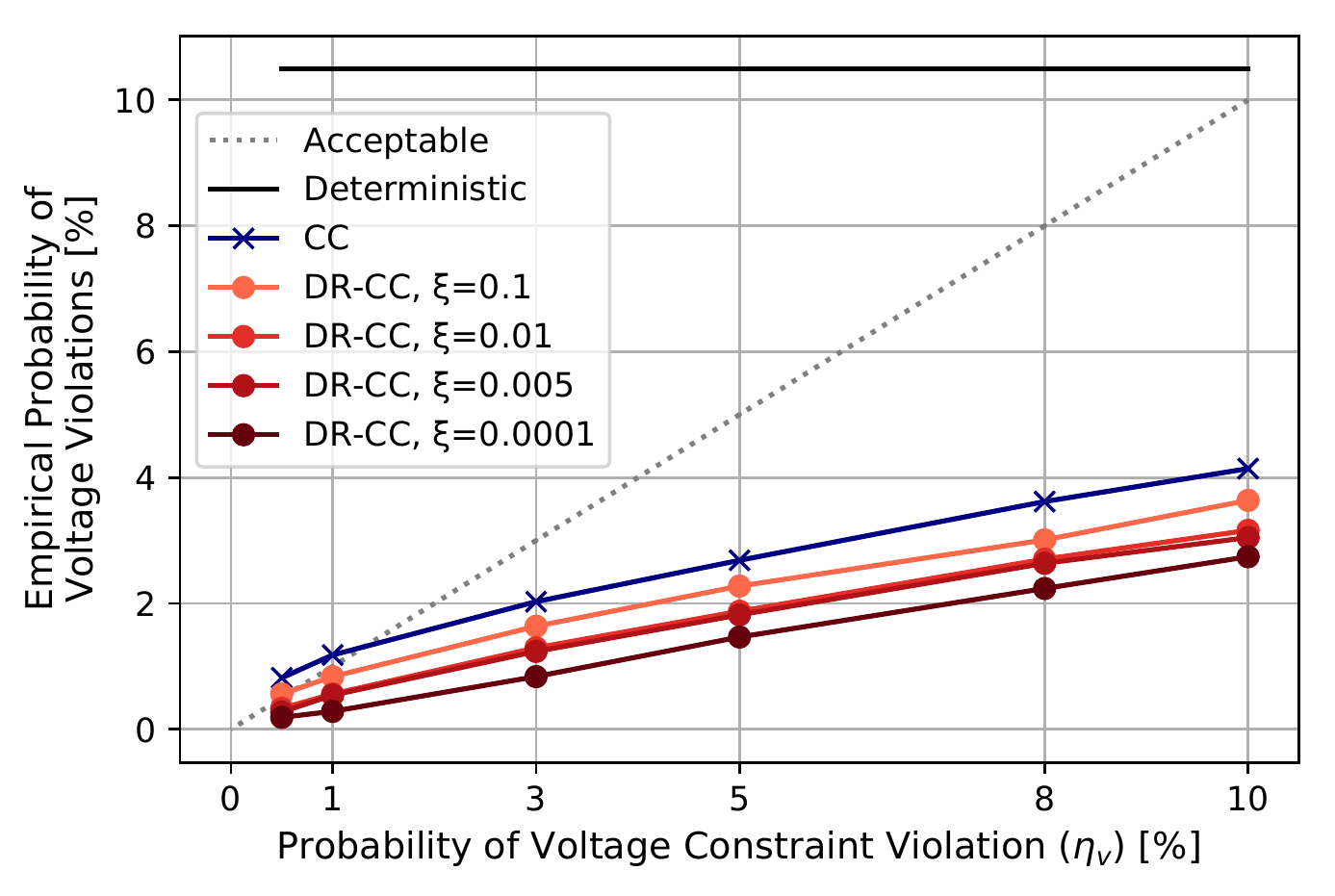}
\caption{Empirical Probability of voltage violation in 750 sample cases for different $\eta_v$ and $\xi$.}
\label{fig:vviolation_to_varconf_and_voltrisk}
\end{figure}

Figure \ref{fig:vviolation_to_varconf_and_voltrisk} presents empirical probabilities of voltage constraint violations (either upper or lower limit) for different $\eta_v$.  As $\eta_v$ increases so does the frequency of observed violations. If $\eta_v>1\%$, the AC-CCOPF  and  DR-AC-CCOPF have lower empirical violations than the postulated value of  $\eta_v$. Note that the DR-AC-CCOPF systematically outpeforms the AC-CCOPF. On the other hand, the deterministic ACOPF systematically underperforms relative to the both CC formulations. As $\xi$ reduces, i.e. the uncertainty set in Fig.~\ref{fig:distribution_visualization}(c) spreads, the DR-AC-CCOPF solution becomes more conservative and returns less violations.

\begin{figure}
\centering
\includegraphics[width=0.8\linewidth]{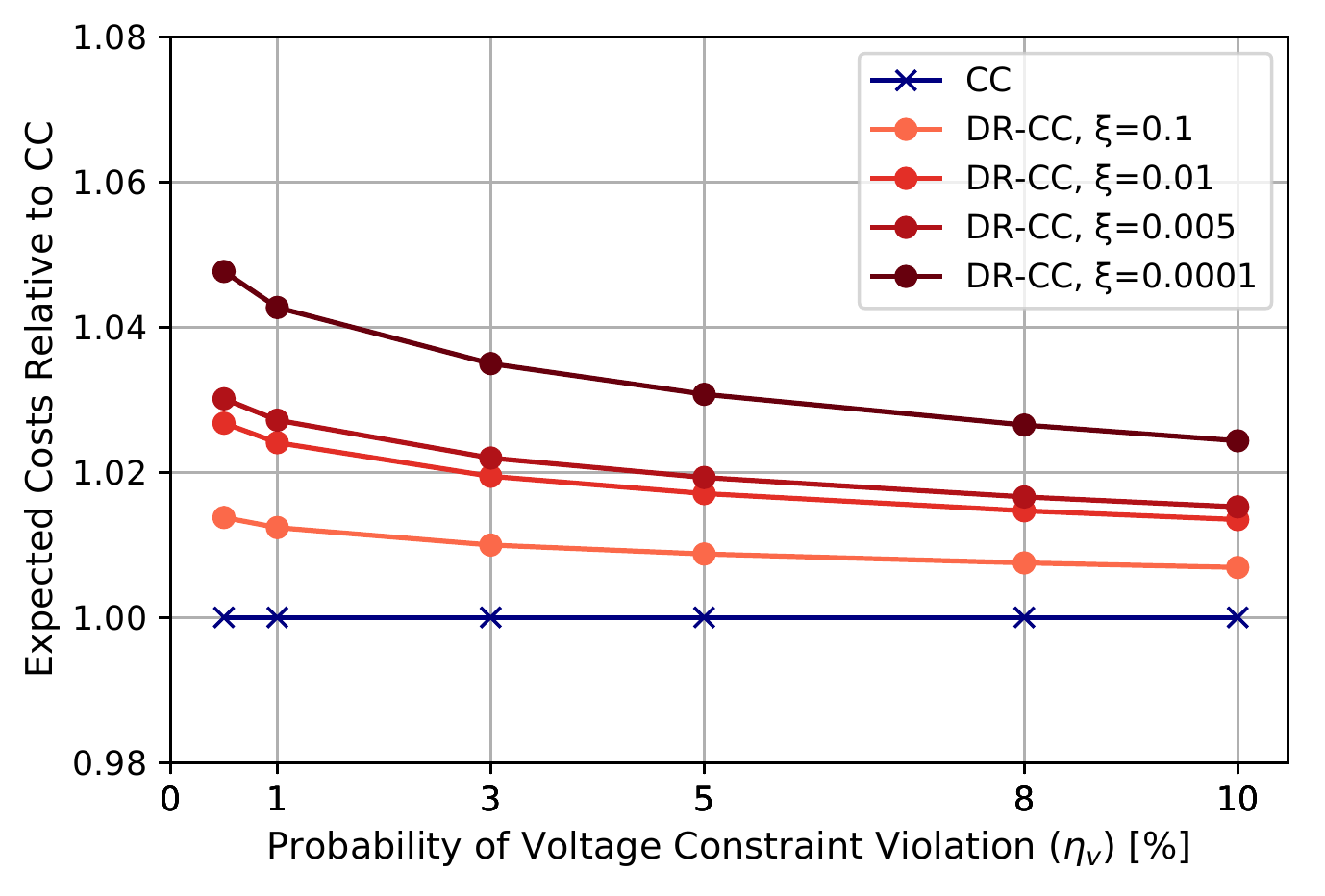}
\caption{Relative expected cost for different $\eta_v$ and $\xi$.}
\label{fig:objective_to_varconf_and_voltrisk}
\end{figure}

Figure \ref{fig:objective_to_varconf_and_voltrisk} compares the expected costs for each OPF formulation normalized by the expected cost of the AC-CCOPF formulation. The natural conservatism of the DR-AC-CCOPF solution, as follows from better compliance with voltage limits in Fig. \ref{fig:vviolation_to_varconf_and_voltrisk},  results in a moderate increase in the expected cost relative to the AC-CCOPF solution. However, the gap between these two formulations narrows as $\eta_v$ increases. Similarly, as the width of the uncertainty set  in Fig. \ref{fig:distribution_visualization}(c) increases, so does the worst-case variance and thus the expected cost as per \eqref{eq:wc_expected_cost}. Our numerical results suggest that the expected cost is more sensitive to changes in $\eta_v$ than to the width of the uncertainty set. The trade-off between the solution feasibility and expected cost is not trivial and DSOs tend to opt a costly, yet more reliable solution, which motivates our out-of-sample performance analysis below.  

\subsection{Out-of-Sample (OOS) Performance}

\begin{figure}[!b]
\centering
\includegraphics[width=0.95\linewidth]{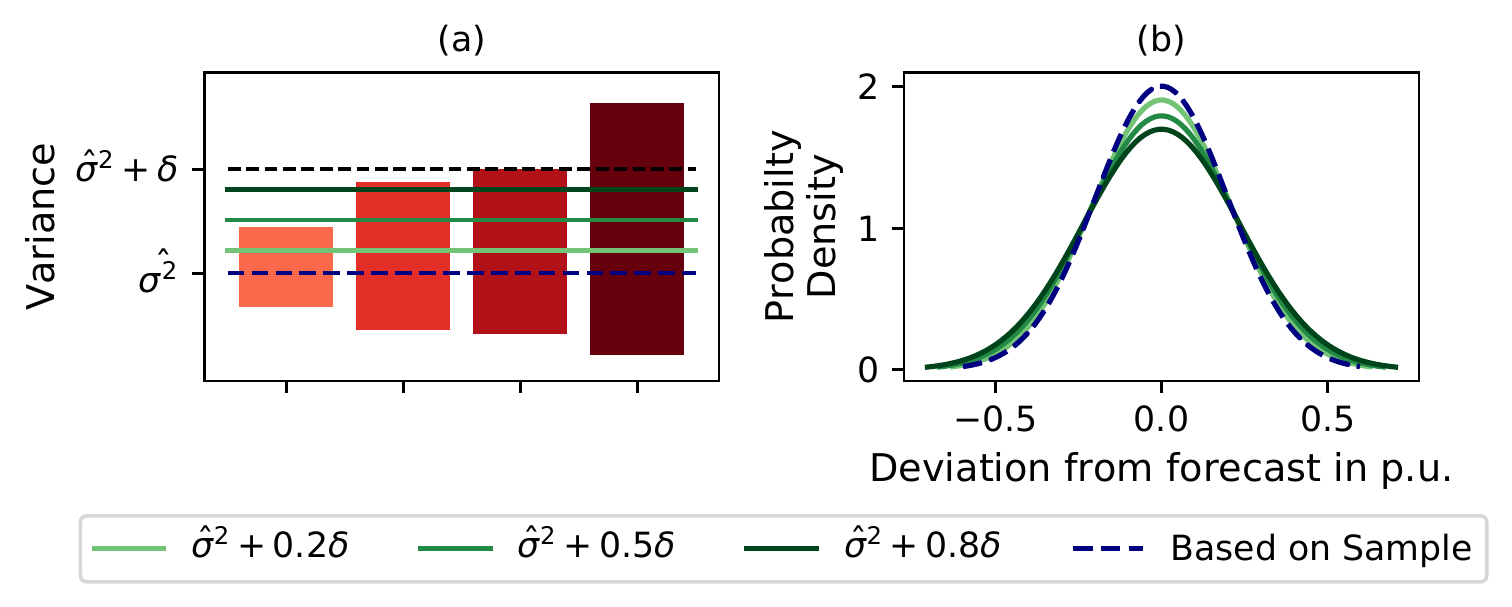}
\caption{Distributions  for out-of-sample testing. \textbf{(a)} Position of out-of-sample variances in the uncertainty intervals. \textbf{(b)}: Resulting error distributions in comparison to the distribution based on sample variance.}
\label{fig:oos_visualization}
\end{figure}

\begin{figure}
\centering
\includegraphics[width=0.8\linewidth]{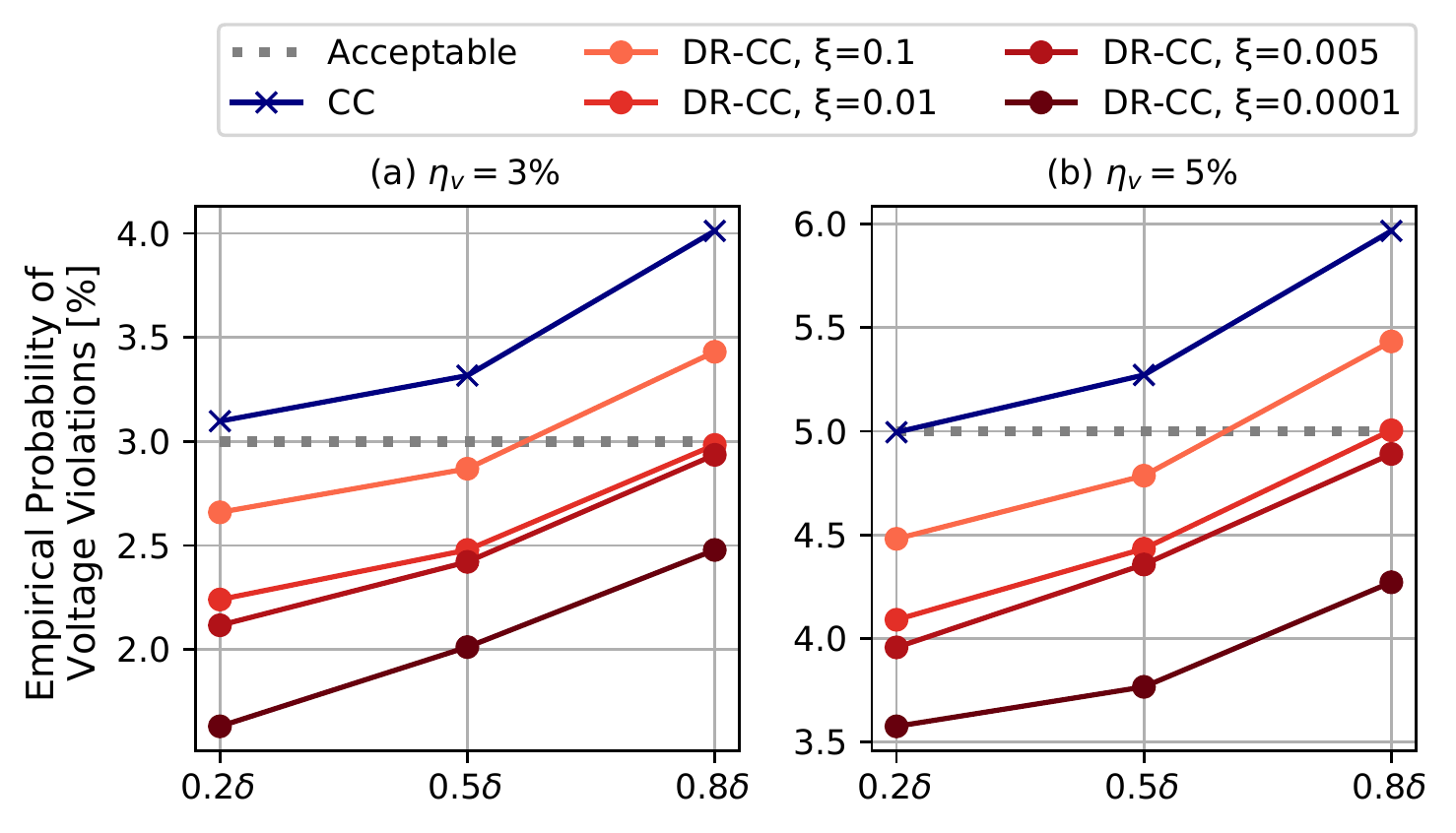}
\caption{Out-of-sample performance: Empirical probability of voltage violations in 750 samples with increasing distance of the true error distribution from the estimated distribution for $\eta_v = 3\%$ (a) and $\eta_v = 5\%$ (b).}
\label{fig:out_of_sample_performance}
\end{figure}

For the following analyses, we generate 750 samples drawn from newly parameterized distributions which are supported by the initial sample data in Figure \ref{fig:distribution_visualization}(b), but have a new value of $\sigma^2$ shifted towards the upper limit of the uncertainty set by parameter $\delta$. Figure~\ref{fig:oos_visualization} shows the three OOS cases and the resulting distributions from which the test samples have been drawn. We use these samples to compare the AC-CCOPF and DR-AC-CCOPF for  $\eta_v = \unit[3]{\%}$ and $\eta_v = \unit[5]{\%}$ as shown in Figure \ref{fig:out_of_sample_performance}. The AC-CCOPF, which is not immunized against distribution ambiguity, does not satisfy the theoretical violation probability limit in all instances, except for the  OOS distribution with $\hat{\sigma}^2 \rightarrow \sigma^2$. On the other hand, the DR-AC-CCOPF holds the theoretical violation probability limit in nearly all cases, except for the two cases related to the smallest distributional uncertainty set. As we can see in Figure \ref{fig:oos_visualization}(a), two of the three OOS cases are outside this set, which explains the empirical violation of the defined $\eta_v$. For other OOS cases, DR-AC-CCOPF meets the requirements posed by the theoretical violation probability limit.

\subsection{Computational Effort and Scalabilty}
The calculations for the case study have been performed on a PC with an Intel Core i5 processor at 2.1GHz with 4GB memory in less than one second each. In order to show the scalabilty of the proposed model, Table \ref{tab:computaion} summarizes computational performance for larger networks based on the IEEE distribution systems data sets%
The results shown in Table~\ref{tab:computaion} use the same values of $\eta_v = 0.05$ and $\xi=0.005$ for all networks. 

\begin{table}[!t]
    \scriptsize
    \centering
    \caption{Computation time for larger systems.}
    \label{tab:computaion}
    \begin{tabular}{r|l}
        \hline
         Case Name & Computing Time (s) \\
        \hline
        15-bus system, \cite{Papavasiliou_Analysis_2017} & $<1$ \\
        IEEE 37-bus system, \cite{schneider2017analytic} & $<1$ \\
        IEEE 123-bus system, \cite{schneider2017analytic} & $1.1$ \\
        IEEE 8500-bus system, \cite{schneider2017analytic} & $24.2$ \\
        \hline
    \end{tabular}
\end{table}

\section{Conclusion}
\label{sec:conclusion}

Leveraging the \emph{LinDistFlow} approximation \cite{Turitsyn_Local_2010,Baran_Optimal_1989} and current research on risk-aware OPF e.g. \cite{Bienstock_Chance_2014}, we have developed and implemented a chance constrained AC OPF formulation for radial distribution systems (AC-CCOPF). To overcome the untenable assumption of perfect knowledge of the underlying probability distributions, the proposed formulation is extended to only rely on historical data. By introducing a distributional uncertainty set and leveraging methods of \emph{distributionally robust optimization}, the formulation is immunized against uncertainty in the probabilistic models of forecast errors obtained from the available observations (DR-AC-CCOPF). The case study reveals that the distributionally robust formulation systemically reduces the empirical probability of voltage violations at a moderate increase in the expected costs. In the conducted out-of-sample performance evaluation, the DR-AC-CCOPF systematically outperforms AC-CCOPF, which is not able to guarantee the theoretical violation probability limit.

\bibliographystyle{IEEEtran}
\bibliography{bibliography_formated.bib}

\end{document}